\documentclass[12pt]{article}
\usepackage{amsmath,amssymb,amsthm,latexsym,amscd}

\title{Ranks for families of all theories \\ of given languages\footnote{
\newline\indent \ \ \ This research was partially
supported by Committee of Science in Education and Science
Ministry of the Republic of Kazakhstan (Grants No. AP05132349,
AP05132546) and Russian Foundation for Basic Researches (Project
No. 17-01-00531-a). } }
\author{Nurlan D.
Markhabatov, Sergey V.
Sudoplatov\footnote{nur\underline{\,\,\,}24.08.93@mail.ru,
sudoplat@math.nsc.ru}}
\date{}
\begin{document}
\maketitle

\begin{abstract}
For families of all theories of arbitrary given languages we
describe ranks and degrees. In particular, we characterize
(non-)totally transcendental families. We apply these
characterizations for the families of all theories of given
languages, with models of given finite or infinite cardinality.

\end{abstract}

The rank \cite{RSrank} for families of theories, similar to Morley
rank, can be considered as a measure for complexity or richness of
these families. Thus increasing the rank by extensions of families
we produce more rich families obtaining families with the infinite
rank that can be considered ``rich enough''.

In the present paper, for families of all theories of arbitrary
given language, we describe ranks and degrees, partially answering
the question in \cite{RSrank}. In particular, we characterize
(non-)totally transcendental families. Thus, we describe rich
families with respect to the rank. Besides, we apply these
characterizations for the families of all theories of given
languages, with models of given finite or infinite cardinality.

\section{Preliminaries}

Throughout we consider families $\mathcal{T}$ of complete
first-order theories of a language $\Sigma=\Sigma(\mathcal{T})$.
For a sentence $\varphi$ we denote by $\mathcal{T}_\varphi$ the
set $\{T\in\mathcal{T}\mid\varphi\in T \}$.

\medskip
{\bf Definition} \cite{at}. Let $\mathcal{T}$ be a family of
theories and $T$ be a theory, $T\notin\mathcal{T}$. The theory $T$
is called {\em $\mathcal{T}$-approximated}, or {\em approximated
by} $\mathcal{T}$, or {\em $\mathcal{T}$-approximable}, or a {\em
pseudo-$\mathcal{T}$-theory}, if for any formula $\varphi\in T$
there is $T'\in\mathcal{T}$ such that $\varphi\in T'$.

If $T$ is $\mathcal{T}$-approximated then $\mathcal{T}$ is called
an {\em approximating family} for $T$, theories $T'\in\mathcal{T}$
are {\em approximations} for $T$, and $T$ is an {\em accumulation
point} for $\mathcal{T}$.

An approximating family $\mathcal{T}$ is called {\em $e$-minimal}
if for any sentence $\varphi\in\Sigma(T)$, $\mathcal{T}_\varphi$
is finite or $\mathcal{T}_{\neg\varphi}$ is finite.

\medskip
It was shown in \cite{at} that any $e$-minimal family
$\mathcal{T}$ has unique accumulation point $T$ with respect to
neighbourhoods $\mathcal{T}_\varphi$, and $\mathcal{T}\cup\{T\}$
is also called {\em $e$-minimal}.

\medskip
Following \cite{RSrank} we define the {\em rank} ${\rm RS}(\cdot)$
for the families of theories, similar to Morley rank
\cite{Morley}, and a hierarchy with respect to these ranks in the
following way.

For the empty family $\mathcal{T}$ we put the rank ${\rm
RS}(\mathcal{T})=-1$, for finite nonempty families $\mathcal{T}$
we put ${\rm RS}(\mathcal{T})=0$, and for infinite families
$\mathcal{T}$~--- ${\rm RS}(\mathcal{T})\geq 1$.

For a family $\mathcal{T}$ and an ordinal $\alpha=\beta+1$ we put
${\rm RS}(\mathcal{T})\geq\alpha$ if there are pairwise
inconsistent $\Sigma(\mathcal{T})$-sentences $\varphi_n$,
$n\in\omega$, such that ${\rm
RS}(\mathcal{T}_{\varphi_n})\geq\beta$, $n\in\omega$.

If $\alpha$ is a limit ordinal then ${\rm
RS}(\mathcal{T})\geq\alpha$ if ${\rm RS}(\mathcal{T})\geq\beta$
for any $\beta<\alpha$.

We set ${\rm RS}(\mathcal{T})=\alpha$ if ${\rm
RS}(\mathcal{T})\geq\alpha$ and ${\rm
RS}(\mathcal{T})\not\geq\alpha+1$.

If ${\rm RS}(\mathcal{T})\geq\alpha$ for any $\alpha$, we put
${\rm RS}(\mathcal{T})=\infty$.

A family $\mathcal{T}$ is called {\em $e$-totally transcendental},
or {\em totally transcendental}, if ${\rm RS}(\mathcal{T})$ is an
ordinal.

\medskip
{\bf Proposition 1.1} \cite{RSrank}. {\em If an infinite family
$\mathcal{T}$ does not have $e$-minimal subfamilies
$\mathcal{T}_\varphi$ then $\mathcal{T}$ is not totally
transcendental.}

\medskip
If $\mathcal{T}$ is totally transcendental, with ${\rm
RS}(\mathcal{T})=\alpha\geq 0$, we define the {\em degree} ${\rm
ds}(\mathcal{T})$ of $\mathcal{T}$ as the maximal number of
pairwise inconsistent sentences $\varphi_i$ such that ${\rm
RS}(\mathcal{T}_{\varphi_i})=\alpha$.

\medskip
Recall the definition of the Cantor--Bendixson rank. It is defined
on the elements of a topological space $X$ by induction: ${\rm
CB}_X(p)\geq 0$ for all $p\in X$; ${\rm CB}_X(p)\geq\alpha$ if and
only if for any $\beta<\alpha$, $p$ is an accumulation point of
the points of ${\rm CB}_X$-rank at least $\beta$. ${\rm
CB}_X(p)=\alpha$ if and only if both ${\rm CB}_X(p)\geq\alpha$ and
${\rm CB}_X(p)\ngeq\alpha+1$ hold; if such an ordinal $\alpha$
does not exist then ${\rm CB}_X(p)=\infty$. Isolated points of $X$
are precisely those having rank $0$, points of rank $1$ are those
which are isolated in the subspace of all non-isolated points, and
so on. For a non-empty $C\subseteq X$  we define ${\rm
CB}_X(C)=\sup\{{\rm CB}_X(p)\mid p\in C\}$; in this way ${\rm
CB}_X(X)$ is defined and ${\rm CB}_X(\{p\})={\rm CB}_X(p)$ holds.
If $X$ is compact and $C$ is closed in $X$ then the sup is
achieved: ${\rm CB}_X(C)$ is the maximum value of ${\rm CB}_X(p)$
for $p\in C$; there are  finitely many points of maximum rank in
$C$ and the number of such points is the \emph{${\rm
CB}_X$-degree} of $C$, denoted by $n_X(C)$.

If $X$ is countable and compact then ${\rm CB}_X(X)$ is a
countable ordinal and every  closed subset has ordinal-valued rank
and finite ${\rm CB}_X$-degree $n_X(X)\in\omega\setminus\{0\}$.

For any ordinal $\alpha$ the set $\{p\in X\mid{\rm
CB}_X(p)\geq\alpha\}$ is called the $\alpha$-th {\em ${\rm
CB}$-derivative} $X_\alpha$ of $X$.

Elements $p\in X$ with ${\rm CB}_X(p)=\infty$ form the {\em
perfect kernel} $X_\infty$ of $X$.

Clearly, $X_\alpha\supseteq X_{\alpha+1}$, $\alpha\in{\rm Ord}$,
and $X_\infty=\bigcap\limits_{\alpha\in{\rm Ord}}X_\alpha$.

It is noticed in \cite{RSrank} that any $e$-totally transcendental
family $\mathcal{T}$ defines a superatomic Boolean algebra
$\mathcal{B}(\mathcal{T})$ with ${\rm RS}(\mathcal{T})={\rm
CB}_{\mathcal{B}(\mathcal{T})}(B(\mathcal{T}))$, ${\rm
ds}(\mathcal{T})=n_{\mathcal{B}(\mathcal{T})}(B(\mathcal{T}))$,
i.e., the pair $({\rm RS}(\mathcal{T}),{\rm ds}(\mathcal{T}))$
consists of Cantor--Bendixson invariants for
$\mathcal{B}(\mathcal{T})$ \cite{BA}.

By the definition for any $e$-totally transcendental family
$\mathcal{T}$ each theory $T\in\mathcal{T}$ obtains the ${\rm
CB}$-rank ${\rm CB}_\mathcal{T}(T)$ starting with
$\mathcal{T}$-isolated points $T_0$, of ${\rm
CB}_\mathcal{T}(T_0)=0$. We will denote the values ${\rm
CB}_\mathcal{T}(T)$ by ${\rm RS}_\mathcal{T}(T)$ as the rank for
the point $T$ in the topological space on $\mathcal{T}$ which is
defined with respect to $\Sigma(\mathcal{T})$-sentences.

\medskip
{\bf Definition} \cite{RSrank}. Let $\alpha$ be an ordinal. A
family $\mathcal{T}$ of rank $\alpha$ is called {\em
$\alpha$-minimal} if for any sentence $\varphi\in\Sigma(T)$, ${\rm
RS}(\mathcal{T}_\varphi)<\alpha$ or ${\rm
RS}(\mathcal{T}_{\neg\varphi})<\alpha$.

\medskip
{\bf Proposition 1.2} \cite{RSrank}. {\em $(1)$ A family
$\mathcal{T}$ is $0$-minimal if and only if $\mathcal{T}$ is a
singleton.

$(2)$ A family $\mathcal{T}$ is $1$-minimal if and only if
$\mathcal{T}$ is $e$-minimal.

$(3)$ For any ordinal $\alpha$ a family $\mathcal{T}$ is
$\alpha$-minimal if and only if ${\rm RS}(\mathcal{T})=\alpha$ and
${\rm ds}(\mathcal{T})=1$. }

\medskip
{\bf Proposition 1.3} \cite{RSrank}. {\em For any family
$\mathcal{T}$, ${\rm RS}(\mathcal{T})=\alpha$, with ${\rm
ds}(\mathcal{T})=n$, if and only if $\mathcal{T}$ is represented
as a disjoint union of subfamilies
$\mathcal{T}_{\varphi_1},\ldots,\mathcal{T}_{\varphi_n}$, for some
pairwise inconsistent sentences $\varphi_1,\ldots,\varphi_n$, such
that each $\mathcal{T}_{\varphi_i}$ is $\alpha$-minimal.}

\section{Ranks for families of theories depending of given languages}

Let $\Sigma$ be a language. If $\Sigma$ is relational we denote by
$\mathcal{T}_\Sigma$ the family of all theories of the language
$\Sigma$. If $\Sigma$ contains functional symbols $f$ then
$\mathcal{T}_\Sigma$ is the family of all theories of the language
$\Sigma'$, which is obtained by replacements of all $n$-ary
symbols $f$ with $(n+1)$-ary predicate symbols $R_f$ interpreted
by $R_f=\{(\bar{a},b)\mid f(\bar{a})=b\}$.

\medskip
{\bf Theorem 2.1.} {\em For any language $\Sigma$ the family
$\mathcal{T}_\Sigma$ is $e$-minimal if and only if
$\Sigma=\emptyset$ or $\Sigma$ consists of one constant symbol.}

\medskip
Proof. If $\Sigma=\emptyset$ or $\Sigma$ consists of one constant
symbol then $\mathcal{T}_\Sigma$ is countable  and consists of
theories $T_n$ with $n$-element models,
$n\in\omega\setminus\{0\}$, and of the theory $T_\infty$ with
infinite models. The theories $T_n$ are finitely axiomatizable by
the sentences witnessing the cardinalities of models and
$T_\infty$ is unique accumulation point for $\mathcal{T}_\Sigma$.
Thus, $\mathcal{T}_\Sigma$ is $e$-minimal.

Now we assume that $\Sigma\ne\emptyset$ and it is not exhausted by
one constant symbol. Below we consider all possible cases.

If $\Sigma$ has a relational symbol $P$ then $\mathcal{T}_\Sigma$
is divided into to infinite definable parts: with empty $P$ and
with nonempty $P$. Therefore, there is a sentence $\varphi$ with
infinite $(\mathcal{T}_\Sigma)_\varphi$ and infinite
$(\mathcal{T}_\Sigma)_{\neg\varphi}$. Hence, $\mathcal{T}_\Sigma$
is not $e$-minimal.

If $\Sigma$ has at least two constant symbols $c_1$ and $c_2$ then
the family $\mathcal{T}_\Sigma$ is divided into two infinite
parts: with $c_1=c_2$ and with $c_1\ne c_2$. It implies that again
$\mathcal{T}_\Sigma$ is not $e$-minimal.

Finally, if $\Sigma$ contains an $n$-ary functional symbol $f$,
$n\geq 1$, then $\mathcal{T}_\Sigma$ is divided into two infinite
parts: with identical $f$ for each element $a$: $f(a,\ldots,a)=a$,
and with $f(a,\ldots,a)\ne a$ for some $a$. It means that again
$\mathcal{T}_\Sigma$ is not $e$-minimal.~$\Box$

\medskip
By Propositions 1.2 and 1.3 each theory $T$ in $e$-minimal
$\mathcal{T}_\Sigma$ has ${\rm RS}_{\mathcal{T}_\Sigma}(T)\leq 1$,
with unique theory having the ${\rm RS}$-rank $1$. Here, following
Theorem 2.1, ${\rm RS}_{\mathcal{T}_\Sigma}(T)= 1$ if and only if
$T$ has infinite models.

\medskip
{\bf Proposition 2.2.} {\em If $\Sigma$ is a language of $0$-ary
predicates then either ${\rm RS}(\mathcal{T}_\Sigma)=1$ with ${\rm
ds}(\mathcal{T}_\Sigma)=2^n$, if $\Sigma$ consists of $n\in\omega$
symbols, or ${\rm RS}(\mathcal{T}_\Sigma)=\infty$, if $\Sigma$ has
infinitely many symbols.}

\medskip
Proof. If $\Sigma$ consists of $n\in\omega$ $0$-ary predicates
$P_1,\ldots,P_n$ then $\mathcal{T}_\Sigma$ has $2^n$ accumulation
points $T_i$ such that each $T_i$ has infinite models and
$(P_1,\ldots,P_n)$ has values
$(\delta_1,\ldots,\delta_n)\in\{0,1\}^n$.

If $\Sigma$ consists of infinitely many $0$-ary predicates $P_i$
then there is an infinite $2$-tree \cite{Spr} formed by
independent values for $P_i$ in $\{0,1\}$, witnessing that there
are no $e$-minimal subfamilies $\mathcal{T}_\varphi$ and producing
${\rm RS}(\mathcal{T}_\Sigma)=\infty$ by Proposition 1.1.~$\Box$

\medskip
Using Proposition 2.2 a totally transcendental family
$\mathcal{T}_\Sigma$, for a language $\Sigma$ of $n$ $0$-ary
predicates, has $2^n$ theories of ${\rm RS}$-rank $1$, each of
which has infinite models.

\medskip
{\bf Proposition 2.3.} {\em If $\Sigma$ is a language of $0$-ary
and unary predicates, with at least one unary symbol $P$, then
either ${\rm RS}(\mathcal{T}_\Sigma)=2^k$ with ${\rm
ds}(\mathcal{T}_\Sigma)=2^m$, if $\Sigma$ consists of $k\in\omega$
unary symbols and $m\in\omega$ $0$-ary predicates, or ${\rm
RS}(\mathcal{T}_\Sigma)=\infty$, if $\Sigma$ has infinitely many
symbols.}

\medskip
Proof. If $\Sigma$ contains $k\in\omega$ unary symbols $P_i$ then
universes can be divided into $2^k$ parts by $P_i$ such that
cardinalities of these parts can vary from $0$ to infinity. So
varying finite cardinalities of the parts we obtain infinitely
many pairwise inconsistent sentences allowing to vary
cardinalities of other parts. Continuing the process for remaining
parts we have $2^n$ steps forming ${\rm
RS}(\mathcal{T}_\Sigma)=2^k$. Having $m\in\omega$ $0$-ary
predicates $Q_j$, sentences witnessing ${\rm
RS}(\mathcal{T}_\Sigma)=2^k$ are implied by $2^m$ pairwise
inconsistent sentences describing values for $Q_j$. Thus, ${\rm
ds}(\mathcal{T}_\Sigma)=2^m$.

If $\Sigma$ contains infinitely many predicate symbols, $0$-ary
and unary, we construct an infinite $2$-tree of sentences formed
by independent values of predicates. Hence, ${\rm
RS}(\mathcal{T}_\Sigma)=\infty$ using Proposition 1.1.~$\Box$

\medskip
In view of Proposition 2.3 there are $2^m$ theories $T$ in
$\mathcal{T}_\Sigma$ having the maximal ${\rm
RS}_{\mathcal{T}_\Sigma}(T)=2^k$. Each such $T$ has only infinite
parts with respect to the predicates $P_i$. Notice also that ${\rm
RS}_{\mathcal{T}_\Sigma}(T)=s\leq 2^k$ if and only if $T$ has
models with exactly $s$ infinite parts.

\medskip
{\bf Proposition 2.4.} {\em If $\Sigma$ is a language of constant
symbols then either ${\rm RS}(\mathcal{T}_\Sigma)=1$ with ${\rm
ds}(\mathcal{T}_\Sigma)=P(n)$, where $P(n)$ is the number for
partitions of $n$-element sets, if $\Sigma$ consists of
$n\in\omega$ symbols, or ${\rm RS}(\mathcal{T}_\Sigma)=\infty$, if
$\Sigma$ has infinitely many symbols.}

\medskip
Proof. If $\Sigma$ consists of constant symbols $c_1,\ldots,c_n$
then we can write in sentences that these constants can
arbitrarily coincide or not coincide. The sentences $(c_i\approx
c_j)^\delta$, $\delta\in\{0,1\}$, define partitions of the set
$C=\{c_1,\ldots,c_n\}$. The number $P(n)$ of these partitions
\cite[Section 5.4]{SO1} defines all possibilities for ${\rm
ds}(\mathcal{T}_\Sigma)$. Since all $\Sigma$-sentences are reduced
to the descriptions $\varphi$ of these partitions as well as to
the descriptions $\psi$ of cardinalities of the sets
$\overline{C}=M\setminus C$, where $M$ are universes of models of
theories in $\mathcal{T}_\Sigma$, we have ${\rm
RS}(\mathcal{T}_\Sigma)=1$, witnessed by $\psi$, and ${\rm
ds}(\mathcal{T}_\Sigma)=P(n)$, witnessed by $\varphi$.

If $\Sigma$ contains infinitely many constant symbols, we
construct an infinite $2$-tree of sentences formed by independent
(in)equalities of constants. Hence, ${\rm
RS}(\mathcal{T}_\Sigma)=\infty$ using Proposition 1.1.~$\Box$

\medskip
By Proposition 2.4, for ${\rm RS}(\mathcal{T}_\Sigma)=1$ there are
$P(n)$ theories $T$ in $\mathcal{T}_\Sigma$ with ${\rm
RS}_{\mathcal{T}_\Sigma}(T)=1$. Each such $T$ is characterized by
existence of infinite models.

\medskip
{\bf Proposition 2.5.} {\em If $\Sigma$ is a language of $0$-ary
and unary predicates, and constant symbols, then either ${\rm
RS}(\mathcal{T}_\Sigma)$ is finite, if $\Sigma$ consists of
finitely many symbols, or ${\rm RS}(\mathcal{T}_\Sigma)=\infty$,
if $\Sigma$ has infinitely many symbols.}

\medskip
Proof. If $\Sigma$ is finite then we can increase ${\rm
RS}(\mathcal{T}_\Sigma)$ till $2^k$ using unary predicates
$P_1,\ldots,P_k$ repeating arguments for Proposition 2.3. The
degree ${\rm ds}(\mathcal{T}_\Sigma)$ is bounded by finitely many
possibilities for values of $0$-ary predicates and for partitions
of constants combining Propositions 2.3 and 2.4.

If $\Sigma$ has infinitely many symbols then it has either
infinitely many $0$-ary predicates, or unary predicates, or
constant symbols. Anyway it is possible to construct an infinite
$2$-tree, as for Propositions 2.3 and 2.4, guaranteeing ${\rm
RS}(\mathcal{T}_\Sigma)=\infty$.~$\Box$

\medskip
As above, ${\rm RS}$-ranks for theories $T$ in a totally
transcendental family $\mathcal{T}_\Sigma$ in Proposition 2.5 are
characterized by the number of infinite $P_i$-parts in models of
$T$.

\medskip
{\bf Proposition 2.6.} {\em If $\Sigma$ is a language containing
an $m$-ary predicate symbol, for $m\geq 2$, or an $n$-ary
functional symbol, for $n\geq 1$, then ${\rm
RS}(\mathcal{T}_\Sigma)=\infty$.}

\medskip
Proof. Using the arguments for the propositions above it suffices
to show that having a binary predicate symbol $Q$ or a unary
functional symbol $f$ it is possible to define infinitely many
independent definable subsets $X_n$, $n\in\omega$, of universes
$M$ for models of theories in $\mathcal{T}_\Sigma$. It is possible
to code these sets $X_n$, even by acyclic directed graphs, by
existence of paths from some elements $a$ without preimages to
elements $b\in X_n$ such that the $(a,b)$-path has the length $n$.
Coding the sets $X_n$ we can form an infinite $2$-tree for
elements in $Y=\bigcup\limits_{n\in\omega}X_n$ such that some
sentences divide $Y$ into continuum many parts by (non)existence
of paths having the lengths $n$. The existence of this $2$-tree
implies that ${\rm RS}(\mathcal{T}_\Sigma)=\infty$ using
Proposition 1.1.~$\Box$

\medskip
{\bf Remark 2.7.} The arguments for Proposition 2.6 allow to
restrict families $\mathcal{T}_\Sigma$ with binary relational
symbols $R$ to the families $\mathcal{T}_{\{R\},{\rm ag}}$ in
graph languages $\{R\}$, of theories of acyclic graphs, and such
that ${\rm RS}(\mathcal{T}_{\{R\},{\rm ag}})=\infty$.

\medskip
Summarizing arguments above we obtain the following theorem.

\medskip
{\bf Theorem 2.8.} {\em For any language $\Sigma$ either ${\rm
RS}(\mathcal{T}_\Sigma)$ is finite, if $\Sigma$ consists of
finitely many $0$-ary and unary predicates, and finitely many
constant symbols, or ${\rm RS}(\mathcal{T}_\Sigma)=\infty$,
otherwise.}

\medskip
Proof. If $\Sigma$ consists of finitely many $0$-ary and unary
predicates, and constant symbols then ${\rm
RS}(\mathcal{T}_\Sigma)$ is finite by Proposition 2.5. Otherwise,
${\rm RS}(\mathcal{T}_\Sigma)=\infty$ by Propositions 2.5 and
2.6.~$\Box$

\section{Application for families of theories depending on cardinalities of models}

The technique for counting of the ranks ${\rm
RS}(\mathcal{T}_\Sigma)$ can be applied for families
$\mathcal{T}_{\Sigma,n}$ of all theories in languages $\Sigma$ and
having $n$-element models, $n\in\omega$, as well as for families
$\mathcal{T}_{\Sigma,\infty}$ of all theories in languages
$\Sigma$ and having infinite models.

Clearly, for any language $\Sigma$,
$\mathcal{T}_\Sigma=\bigcup\limits_{n\in\omega}\mathcal{T}_{\Sigma,n}\cup\mathcal{T}_{\Sigma,\infty}$.
Therefore, by monotony of ${\rm RS}$, we have for any
$n\in\omega$:
\begin{equation}\label{feq1}{\rm
RS}(\mathcal{T}_{\Sigma,n})\leq {\rm
RS}(\mathcal{T}_\Sigma),\end{equation}
\begin{equation}\label{feq2}{\rm
RS}(\mathcal{T}_{\Sigma,\infty})\leq {\rm
RS}(\mathcal{T}_\Sigma).\end{equation} Using (\ref{feq1}) and
(\ref{feq2}), the following theorems and their arguments allow to
count the ranks ${\rm RS}(\mathcal{T}_{\Sigma,n})$ and ${\rm
RS}(\mathcal{T}_{\Sigma,\infty})$ depending on $\Sigma$.

\medskip
{\bf Theorem 3.1.} {\em For any language $\Sigma$ either ${\rm
RS}(\mathcal{T}_{\Sigma,n})=0$, if $\Sigma$ is finite or $n=1$ and
$\Sigma$ has finitely many predicate symbols, or ${\rm
RS}(\mathcal{T}_{\Sigma,n})=\infty$, otherwise.}

\medskip
Proof. If $\Sigma$ is finite then $\mathcal{T}_{\Sigma,n}$ is
finite for any $n\in\omega$, since there are finitely many
isomorphism types for $n$-element structures in the language
$\Sigma$. If $n=1$ and $\Sigma$ has finitely many predicate
symbols then again there are finitely many isomorphism types for
$1$-element structures $\langle A;\Sigma\rangle$, since there are
finitely many possibilities for distributions of empty predicates,
all nonempty predicates are complete, all constants has same
interpretations, and all functions are identical.

If $\Sigma$ has infinitely many predicate symbols $P_i$, we can
form an infinite $2$-tree of sentences allowing $P_i$
independently be empty or complete. If $\Sigma$ has infinitely
many constant symbols $c_i$, then, for $n\geq 2$ and $c_0\ne c_1$,
we again can form an infinite $2$-tree of sentences allowing $c_i$
independently be equal to $c_0$ or $c_1$. Finally, if $\Sigma$ has
infinitely many functional symbols $f_i$, then, for $n\geq 2$, we
can form an infinite $2$-tree of sentences allowing $f_i$ be
(non)identical. Each possibility above immediately implies ${\rm
RS}(\mathcal{T}_{\Sigma,n})=\infty$.~$\Box$

\medskip
{\bf Theorem 3.2.} {\em For any language $\Sigma$ either ${\rm
RS}(\mathcal{T}_{\Sigma,\infty})$ is finite, if $\Sigma$ is finite
and without predicate symbols of arities $m\geq 2$ as well as
without functional symbols of arities $n\geq 1$, or ${\rm
RS}(\mathcal{T}_{\Sigma,n})=\infty$, otherwise.}

\medskip
Proof. Let $\Sigma$ be finite and without predicate symbols of
arities $m\geq 2$ as well as without functional symbols of arities
$n\geq 1$, i.e., $\Sigma$ contains only finitely many $0$-ary and
unary predicate symbols as well as finitely many constant symbols.
Then applying Propositions 2.2--2.5 and the inequality
(\ref{feq2}) we have ${\rm
RS}(\mathcal{T}_{\Sigma,\infty})<\omega$.

If $\Sigma$ has predicate symbols of arities $m\geq 2$ or
functional symbols of arities $n\geq 1$ then ${\rm
RS}(\mathcal{T}_{\Sigma,n})=\infty$ repeating arguments for
Proposition 2.6 and constructing a $2$-tree of sentences.

If $\Sigma$ is infinite then by the previous case it suffices to
consider languages with either infinitely many $0$-ary predicates,
or infinitely many unary predicates, or infinitely many constants.
In these cases we repeat arguments for Propositions 2.2--2.5 and
construct $2$-trees of sentences guaranteeing ${\rm
RS}(\mathcal{T}_{\Sigma,n})=\infty$.~$\Box$

\medskip
Notice that, similar the remark after Proposition 2.5, ${\rm
RS}$-ranks for theories $T$ in a totally transcendental family
$\mathcal{T}_{\Sigma,\infty}$ are characterized by numbers of
infinite parts, in models of $T$, with respect to unary
predicates.


\end{document}